\documentclass[11pt,leqno]{article}
\usepackage{graphicx}
\usepackage{amsfonts, amsthm}
\begin{document}
\newtheorem{theo}{Theorem}
\newtheorem{exam}{Example}
\newtheorem{coro}{Corollary}
\newtheorem{defi}{Definition}
\newtheorem{prob}{Problem}
\newtheorem{lemm}{Lemma}
\newtheorem{prop}{Proposition}
\newtheorem{rem}{Remark}
\newcommand\Pn[1]{\mathbb{P}^{#1}}  
\newcommand\Hi[1]{\mathbb{P}^{#1}_\infty} 
\def\Z{\mathbb{Z}} 
\def\Q{\mathbb{Q}} 
\def\C{\mathbb{C}}   
\def\deg{{\rm deg}}  
\def\pt{\mathbb{C}[t]} 
\def\N{\mathbb{N}}
\def\H{\mathbb{H}}
\def\W{{\cal W}}     
\def\F{{\cal F}}     
\def\HH{{\cal H}}    
\begin{center}
{\LARGE\bf
Mixed Hodge structure of affine hypersurfaces}
\footnote{ Math. classification: 14C30,
32S35 
} \\
\vspace{.25in} {\large {\sc Hossein Movasati}} \\
\end{center}
\begin{abstract}
In this article 
we give an algorithm which produces a basis
of the $n$-th de Rham cohomology of the 
affine smooth hypersurface $f^{-1}(t)$
compatible with the mixed Hodge structure,
 where $f$ is a polynomial in $n+1$ variables and
satisfies a certain regularity
condition at infinity (and hence has isolated
singularities). 
As an application we show  that the notion of a
Hodge cycle in regular fibers of $f$ is given 
in terms of the vanishing of integrals of
certain polynomial $n$-forms in $\C^{n+1}$ over topological
$n$-cycles on the fibers of
$f$. Since the $n$-th
homology of a regular fiber is generated by vanishing
cycles, this leads us to study Abelian integrals over them.
Our result generalizes and uses the arguments of  J.
Steenbrink in \cite{st77} for quasi-homogeneous polynomials.
\end{abstract}
\setcounter{section}{-1}
\section{Introduction}
\def\O{{\cal O}}
To study the monodromy and some numerical invariants of a
local holomorphic function  $f:(\C^{n+1},0)\rightarrow (\C,0)$ with an 
isolated critical point at $0$,
E. Brieskorn in \cite{bri} introduced a $\O_{\C,0}$ module $H'$ and
the notion of Gauss-Manin connection on $H'$. Later
J. Steenbrink \cite{st76},
inspired by  P. Deligne's theory of mixed Hodge structures
(see \cite{de71} and two others with
the same title) on  algebraic varieties defined over complex numbers
and W. Schmid's limit Hodge structure (see \cite{sc73})
associated to a fibration with projective manifolds,
introduced the notion of the  mixed Hodge structure for a germ of
a singularity $f$. A different construction of such a  mixed Hodge structure
was also given by A. Varchenko in \cite{var} using the asymptotic of integrals of holomorphic forms 
over vanishing cycles.  


In the case of a polynomial $f$ in $\C^{n+1}$, on the  one hand
the $n$-th cohomology of a regular fiber carries  Deligne's mixed Hodge
structure and on the other hand we have the Brieskorn module $H'$ of $f$
which contains the information of the $n$-th de Rham cohomology
of regular fibers. 
The variation of mixed Hodge structures in such situations is studied by 
Steenbrink and Zucker (see \cite{stzu}). In this article we define two filtrations
on $H'$ based on the mixed Hodge structure of the regular fibers of $f$. 
At the beginning my purpose was
to find explicit descriptions
of  arithmetic properties of Hodge cycles
for hypersurfaces in projective spaces. Such descriptions for CM-Abelian
varieties are well-known but in the case of hypersurfaces
we have only descriptions for Fermat varieties
(see \cite{sh79}).  As an application we will see that
it is possible to  write down the property of being a Hodge cycle in terms of
the vanishing of certain integrals 
over cycles generated by vanishing cycles. 
Such integrals also
appear in the context of
holomorphic foliations/differential equations (see \cite{mov0, mov}
and the references there).
The first advantage 
of this approach is that
we can write some consequences of the Hodge conjecture in terms
of periods
(see the example at the end of this Introduction and \S \ref{examples}).
We explain below the  results in this article.

Let $\alpha=(\alpha_1,\alpha_2,\ldots,\alpha_{n+1})\in\N^{n+1}$  
and assume that
the greatest common divisor of all the $\alpha_i$'s is one.
We consider $\C[x]:=\C[x_1,x_2,\ldots,x_{n+1}]$ as a graded
algebra with $deg(x_i)=\alpha_i$.  A polynomial $f\in\C[x]$ is called
 a quasi-homogeneous polynomial of degree $d$ with respect to the grading
$\alpha$ if
$f$ is a linear combination of monomials of the type
$x^\beta:=x_1^{\beta_1}x_2^{\beta_2}\cdots
x_{n+1}^{\beta_{n+1}},\alpha.\beta:=\sum_{i=1}^{n+1}\alpha_i\beta_i=d$. For
an arbitrary polynomial $f\in\C[x]$ one can write in a unique way
$f=\sum_{i=0}^df_i,\ f_d\not=0$,
where $f_i$ is a quasi-homogeneous polynomial of degree
$i$. The number $d$  is called the degree of
$f$. Set
${\mathrm Sing}(f):=\{\frac{\partial f}{\partial x_i}=0,
\ i=1,2,\ldots,n+1\}$.

Let us be given a polynomial $f\in\C[x]$.
We assume that $f$ is a (weighted) strongly tame polynomial. In this article this means
that  there exist natural numbers $\alpha_1,\alpha_2,\ldots,\alpha_{n+1}\in\N$ such that
${\mathrm Sing}(g)=\{0\}$, where $g=f_d$ is the last quasi-homogeneous piece
of $f$ in the graded algebra $\C[x],\ {\mathrm deg}(x_i)=\alpha_i$.
Looking at $f$ as a rational function in the weighted projective
space (see \S \ref{wps}) we will see that the strongly tameness condition on $f$
implies that the polynomial $f$ has isolated singularities,
i.e. ${\mathrm Sing}(f)$ is a discrete set in $\C^{n+1}$.

We choose a basis
$x^I=\{x^\beta\mid\beta\in I\}$ of monomials for the Milnor $\C$-vector space
$V_g:=\C[x]/<\frac{\partial g}{\partial
  x_i}\mid i=1,2,\ldots,n+1>$ and define
\begin{equation}
\label{21may2004}
w_i:=\frac{\alpha_i}{d},\ 1\leq i\leq n+1,\ \eta:=(\sum_{i=1}^{n+1}(-1)^{i-1}
w_i x_i\widehat {dx_i}), A_\beta:=
\sum_{i=1}^{n+1}(\beta_i+1)w_i,\
\eta_\beta:=x^\beta \eta, \ \beta\in I
\end{equation}
where $\widehat{dx_i}=dx_1\wedge\cdots
\wedge dx_{i-1}\wedge dx_{i+1}\wedge\cdots\wedge dx_{n+1}$. It turns out
that $x^I$ is also a basis of $V_f$ and so $f$ and $g$ have the same Milnor 
numbers (see \S \ref{basis}). We denote  it by  $\mu$ . 
We denote by $C=\{c_1,c_2,\ldots,c_r\}\subset\C$
the set of critical values of $f$ and $L_c:=f^{-1}(c),\ c\in\C$.
The strongly tameness condition on $f$ implies that the fibers 
$L_c,c\in \C\backslash C$ are connected and the function $f$ is a 
$C^\infty$ fiber bundle on $\C\backslash C$ (see \S \ref{wps}).
Let $\Omega^i, i=1,2,\ldots, n+1$ be the set of polynomial differential
$i$-forms in $\C^{n+1}$.
The Brieskorn module
$$
H'=H'_f:=\frac
{\Omega^n}{df\wedge \Omega^{n-1}+d\Omega^{n-1}}
$$
of $f$ is a $\C[t]$-module
in a natural way: $t.[\omega]=[f\omega],\ [\omega]\in H'$. If there is no danger of confusion 
we will not write the brackets. A direct 
generalization of the topological argument in \cite{ga98} in the case $n=1$ 
implies that $H'$ is freely generated by  the forms 
$\eta_\beta,\ \beta\in I$ (see Proposition \ref{28.6.04} for an algebraic 
proof).

Using vanishing theorems and the Atiyah-Hodge-Grothendieck theorem on
the de Rham cohomology of affine varieties (see \cite{grot}), we see that 
$H'$ restricted to each regular fiber $L_c,\ c\in\C\backslash C$ is 
isomorphic to the $n$-th de Rham  cohomology  of 
$L_c$ with complex coefficients. The Gauss-Manin connection associated to 
the fibration $f$ on $H'$ turns out to be a map 
$$
\nabla: H'\rightarrow H'_C 
$$
satisfying the Leibniz rule, where for a set $\tilde C\subset\C$ by 
$H'_{\tilde C}$ we mean the localization of $H'$ on the
multiplicative subgroup of $\C[t]$ generated by $t-c,\ c\in \tilde C$ 
(see \S \ref{brieskorn}). 
Using the Leibniz rule one can extend $\nabla$ to a function from 
$H'_{\tilde C}$ 
to itself.  Here $\tilde C$ is any subset of $\C$ containing $C$.

The mixed Hodge structure $(W_\bullet, F^\bullet)$ of $H^n(L_c,\C)$ 
is defined  by Deligne in \cite{de71} using the hypercohomology interpretation 
of  the cohomology of $L_c$  and the sheaf of meromorphic
forms with logarithmic poles. It is natural to define  a  double filtration $(W_\bullet, F^\bullet)$
in $H'$ as follows: $W_mH',\ m\in\Z$ (resp. $F^kH',\ k\in\Z$) 
consists of elements $\psi \in H'$ such that the restriction of 
$\psi$ on all $L_c$'s except a finite number of them belongs to 
$W_mH^n(L_c,\C)$ (resp. $F^kH^n(L_c,\C)$).
In connection to the work of Steenbrink and Zucker, we mention that
on $\C\backslash C$ we have the variation of mixed Hodge structures 
$H^n(L_t,\C), \ t\in \C\backslash C$. 
The Brieskorn module $H'$ for a strongly tame polynomial gives a free extension to 
$\C$ of the underlying free $\O_{\C\backslash C}$-module. 
Here we identify coherent sheaves on $\C$ with finite modules
over $\C[t]$ by taking the global sections. 
Therefore, the mentioned filtrations of $H'$ in this
text are given by the maximal extensions as $\C[t]$-submodules of the
Brieskorn module. 
Since in our situation $H'$ is freely generated of finite rank,  they 
are also freely generated sub-modules of $H'$. 
Their rank is equal to to the dimensions of the mixed Hodge structure 
of a regular  fiber of $f$. Note that we do not know yet whether 
$Gr^{k}_F Gr^W_m H',\ k\in\Z,\ m=n,n+1$ 
are freely  generated $\C[t]$-modules. 
In the same way we define $(W_\bullet,F^\bullet)$ of 
the localization of $H'$ over multiplicative subgroups of $\C[t]$.
In this article we prove that:
\begin{theo}
\label{main}
Let $b\in\C\backslash C$ be a regular value of $f\in\C[x]$. If $f$ 
is a (weighted) strongly tame polynomial then $Gr_{m}^WH'=0$ for 
$m\not =n,n+1$ and there exist 
a map $\beta\in I\rightarrow d_\beta\in\N\cup\{0\}$  and $C\subset \tilde 
C\subset \C$ such that $b\not\in \tilde C$ and
\begin{equation}
\label{naucia}
\nabla^{k} \eta_\beta,\ \beta\in I,\ A_\beta=k
\end{equation}
form a basis of $Gr^{n+1-k}_FGr_{n+1}^W H'_{\tilde C}$
and the forms
\begin{equation}
\label{vomiting}
\nabla^{k}\eta_\beta,
A_\beta+\frac{1}{d}\leq k
\leq A_\beta +\frac{d_\beta}{d}
\end{equation}
form a basis of $Gr^{n+1-k}_F Gr^W_n H'_{\tilde C}$, where 
$\nabla^k=\nabla\circ\nabla\circ\cdots\circ\nabla$ $k$ times.
\end{theo}
The numbers $d_\beta$  and the set $\tilde C$ are calculated from a monomial
basis of the Jacobian of the homogenization of $f-b$ 
(see Lemma \ref{10:35}) and hence they
are  not unique and  may depend on the choice of $b$.
For a generic $b$ one can put $d_\beta=d-1$ but this is not the case for
all $b$'s (see Example \ref{5.7.03}).
In \cite{mov04} the equality $\tilde C=C$  is shown for many examples of 
$f$ in two variables and for a suitable choice
of $d_\beta$'s. For those examples a similar theorem is proved
as above for  the Brieskorn module rather than its localization. This has many applications in the theory of Abelian integrals
in differential equations (see \cite{ga02, mov}).   
When $f=g$ is a quasi-homogeneous polynomial
of degree $d$ with an isolated singularity at $0\in\C^{n+1}$ our result
can be obtained from J. Steenbrink's result in \cite{st77} using the residue
theory adapted
to Brieskorn modules (see Lemma \ref{hayhay}, \S \ref{residue}).
In this case $C=\{0\}$ and
any two regular fibers are biholomorphic. We have $d_\beta=d-1,\ \forall \beta\in I$, $\nabla(\eta_\beta)=
\frac{A_\beta}{t}\eta_\beta$ and so 
$\nabla^{k}\eta_\beta=\frac{A_\beta(A_\beta-1)\cdots (A_\beta-k+1)}{t^k}
\eta_\beta$. In this case we get the following stronger statement:  
$\eta_\beta, A_\beta=k\in\N$ form a basis 
of $Gr_F^{n+1-k}Gr^W_{n+1}H'$ and $\eta_\beta, A_\beta\not\in \N, 
-[-A_\beta]=k$ form
a basis of $Gr_F^{n+1-k}Gr^W_{n}H'$. 

One may look at the fibration $f=t$ as an affine variety $X$
defined over the  function field $\C(t)$ and interpret Theorem \ref{main} 
as the existence of mixed Hodge structure on the de Rham cohomology of $X$ 
(see \cite{grot} and also \cite{pin}).
However, we note that the Brieskorn module is something
finer than the de Rham cohomology of $X$; for instance if we do not have 
the tameness 
property $H'$ may not
be finitely generated.

One of the initial motivations for me to get  theorems
like Theorem \ref{main} was in obtaining the property of
being a
Hodge cycle in terms of the  vanishing of explicit integrals of 
polynomial $n$-forms in
$\C^{n+1}$. In the case $n$ even,
a cycle in $H_n(L_c,\Z),\ c\not\in C$ is called a Hodge cycle if
its image in $H_n(\widehat{L_c},\Z)$ is a Hodge cycle,
where $\widehat{L_c}$ is the  smooth compactification
of $L_c$. Since the mixed Hodge structure on $H^n(L_c,\C)$ is independent
of the compactification and the map
$i: H^n(\widehat L_c,\C)\rightarrow H^n(L_c,\C)$ induced by the
inclusion $L_c\subset \widehat{L_c}$ is a weight zero
morphism of mixed Hodge structures, this definition  does not depend
on the compactification of $L_c$. Moreover, $Gr_n^WH^n(L_c,\C)$ in the case 
$\alpha_1=\cdots=\alpha_{n+1}=1$ coincides
with the primitive cohomology of the canonical compactification and hence, 
we capture all the Hodge cycles contained in the primitive cohomology (via 
Poincar\'e duality). 


\begin{coro}
\label{261103}
In the situation of Theorem \ref{main}, 
a cycle $\delta_c\in H_n(L_c,\Z), c\in \C\backslash \tilde C$ is Hodge 
if and only if
$$
(\frac{\partial^k}{\partial t^k}\int_{\delta_t}\eta_\beta)\mid_{t=c}=0,\ \forall 
(\beta,k) \in I_h,
$$
where $I_h=\{(\beta,k)\in I\times\Z \mid  A_\beta+\frac{1}{d}\leq k
\leq A_\beta +\frac{d_\beta}{d},\ A_\beta\not \in\N,\ k\leq \frac{n}{2} \}$ 
and $\{\delta_t\}_{t\in(\C,c)}$ is a continuous family of cycles
in the fibers of $f$ which is obtained by using a local topological
trivialization of $f$ around $c$.
\end{coro}
Note that 
$$(\frac{\partial^k}{\partial t^k}\int_{\delta_t}\eta_\beta)\mid_{t=c}=
\int_{\delta_c} \nabla^k\eta_\beta=\sum_{\beta\in I}p_{\beta,k}(c)\int_{\delta_c}
\eta_\beta
$$ 
where $\nabla^k\eta_\beta=\sum_{\beta\in I}p_{\beta,k}\eta_\beta,\ p_{\beta,k}\in\C[t]_C$, and the forms $\nabla^k\eta_\beta,(\beta,k)\in I_h$ form a basis 
of $F^{\frac{n}{2}+1}H'_{\tilde C}\cap W_{n}H'_{\tilde C}$.
Also $H_n(L_c,\Z)$ is generated by a distinguished set
of vanishing cycles
(see \cite{dine, arn}) and one may be
interested in constructing such  a distinguished set of vanishing cycles, 
try to carry out explicit integration and hence obtain more explicit
descriptions of Hodge cycles.
For an $\omega \in H'$ the function $h(t)=\int_{\delta_t}\omega$ 
extends  to a multivalued function on $\C\backslash C$ and satisfies
a Picard-Fuchs equation with possible poles at $C$.   
For a quasi-homogeneous polynomial $f=g$ the Picard-Fuchs equation associated 
to $\eta_\beta$ is $t\frac{\partial h}{\partial t}-A_\beta.h=0$.
For the example $f=x_1^3+x_2^3+\cdots+x_5^3-x_1-x_2$ which has 
a non-trivial variation of Hodge structures using 
{\sc Singular}  (see \S \ref{examples}) we get the following fact:  
For all $c\in\C-\{\pm \frac{4}{3\sqrt{3}}, 0\}$ a
cycle $\delta\in H_4(L_c,\Z)$ is Hodge if and only if 
{\tiny
$$
(972c^2-192)\int_\delta x_1x_2\eta+
(-405c^3-48c)\int_\delta x_2\eta
+(-405c^3-48c)\int_\delta x_1\eta
+(243c^4-36c^2+64)\int_\delta \eta=0
$$
}
Since the Hodge conjecture is proved for cubic hypersurfaces of dimension $4$ 
by  C. Clemens, J. P. Murre and S. Zucker (see \cite{zu77}), we 
conclude that all the integrals appearing in the above equality 
divided by $\pi^2$ are algebraic numbers (see \S \ref{examples}). 
In \cite{stho} we have used this idea to compute the values of  the
Gauss hypergeometric series at certain algebraic points. 

Let us explain the structure of the article. 
In \S \ref{wps} we
recall some terminology on weighted projective spaces.  
In \S \ref{hodge} we explain the idea that to
be able to give descriptions of Hodge cycles in terms of integrals  one must
consider them with support in affine varieties and then
use Theorem \ref{main} and get the property of being a Hodge
cycle in terms of the  vanishing of integrals.
In \S \ref{brieskorn} we introduce two Brieskorn modules $H'$ and
$H''$ associated to a polynomial $f$ and the notion of Gauss-Manin
connection on them. The version of Gauss-Manin connection we use
here comes from the context of differential equations (see
\cite{mov}) and the main point about it is that we can iterate it.
In \S \ref{residue} we see how the iteration of an element
$\omega$ of $H''$ by the Gauss-Manin connection is related to the
residue of $\omega/(f-c)^k,k\in \N$ in the regular fiber $L_c$ of $f$.
\S \ref{grst} is dedicated to a
generalization of a theorem of Griffiths (see \cite{gr69}) to
weighted projective spaces by J. Steenbrink (see \cite{st77}). 
The main point in this section is Theorem \ref{1.2.04}.  
What is new is an explicit basis of the underlying cohomology.
In \S \ref{basis} we prove Theorem \ref{main}.
\S \ref{examples}  is dedicated to some examples. 

When the first draft of this paper was written M. Schulze told me about his 
article \cite{sch} 
in which he gives an algorithm to calculate a good 
$\C[\nabla^{-1}]$-basis  of the Brieskorn module for strongly tame polynomials. 
For the moment the only thing which I can say is that in the case 
$f=g$ the set $\{\eta_\beta,\beta\in I\}$ is also a good basis of
the $\C[\nabla^{-1}]$-module $H'$ because 
$t\eta_\beta=(A_\beta+1)\nabla^{-1}\eta_\beta, \ \forall \beta\in I$. 
In particular, any strongly tame polynomial in the sense of this article has the
same monodromy at infinity as $g$ and so the spectrum of $H'$ is equal to
$\{A_\beta+1,\beta\in I\}$. The generalization of the result of this article
for a tame polynomial in the sense of \cite{sab}  or a Lefschetz pencil (see \cite{momo}) 
would be a nice challenge. Note that the pair $(W_\bullet,F^\bullet)$ 
defined on $H'$ is different form
the mixed Hodge structure constructed  in \cite{sab} (using also the Gauss-Manin connection), 
in which the weight filtration is due to  the monodromy at infinity. 

One can compute $\nabla$, $d_\beta$'s and calculate every element of $H'$ 
as a $\C[t]$-linear combination of $\eta_\beta$'s. 
These are done in the library {\tt foliation.lib} written in {\sc Singular} 
and is explained in \cite{mov04}. 
This paper is devoted to the applications
of the existence of such a basis in differential equations. 
Also for many examples it is shown that by modification of 
Theorem \ref{main} one can get a basis of $H'$ 
compatible with $W_\bullet H'$ and $F^\bullet H'$.

{\bf Acknowledgment:} I learned Hodge theory when I was at the
Max-Planck Institute for Mathematics in Bonn and in this direction
S. Archava helped me a lot. Here I would like to thank him and
the Institute.
When I had the rough idea of the results of this article in my mind, I visited
Kaiserslautern, where  Gert-Martin Greuel drew my attention to
the works of J. Steenbrink. Here I would like to thank him and 
the {\sc Singular} team. 
I would like to thank the Mathematics Department
of the University of G\"ottingen, where the main result of this
article was obtained, for hospitality and financial support.
My thanks also go to T. E. Venkata Balaji for useful conversations
in Algebraic Geometry.
\section{Weighted projective spaces}
\label{wps}
In this section we recall some terminology on weighted projective spaces.
For more information on weighted projective spaces the
reader is referred to \cite{do82, st77}.

Let $n$ be a natural number and
$\alpha:=(\alpha_1,\alpha_2,\ldots,\alpha_{n+1})$ be a vector of natural
numbers
whose greatest common divisor
is one.
The multiplicative group $\C^*$ acts on
$\C^{n+1}$ in the following way:
$$
(X_1,X_2,\ldots,X_{n+1}) \rightarrow
(\lambda^{\alpha_1}X_1,\lambda^{\alpha_2}X_2,
\ldots,\lambda^{\alpha_{n+1}}X_{n+1}),\ \lambda\in\C^*
$$
We also denote the above map by $\lambda$. The quotient space
$$
\Pn {\alpha}:= \C^{n+1}/\C^*
$$
is
called the projective space of weight $\alpha$.
If $\alpha_1=\alpha_2=\cdots=\alpha_{n+1}=1$ then $\Pn \alpha$
is the usual projective
space $\Pn {n}$ (Since $n$ is a natural number,
$\Pn n$ will not mean a zero dimensional weighted projective space).
One can give another interpretation of $\Pn \alpha$ as
follow: Let $G_{\alpha_i}:=\{ e^{\frac{2\pi\sqrt{-1}m }{\alpha_i}}
\mid m\in\Z\}$.
The group $\Pi_{i=1}^{n+1}
G_{\alpha_i}$ acts discretely on the usual projective space $\Pn n$ as
follows:
$$
(\epsilon_1,\epsilon_2,\ldots,\epsilon_{n+1}), [X_1:X_2:\cdots:X_{n+1}]
\rightarrow
[\epsilon_1X_1:\epsilon_2X_2:\cdots:\epsilon_{n+1}X_{n+1}]
$$
The quotient space $\Pn n/ \Pi_{i=1}^{n+1} G_{\alpha_i}$ is
canonically isomorphic to $\Pn \alpha$. This canonical isomorphism is
given by $$
[X_1:X_2:\cdots:X_{n+1}]\in \Pn n/ \Pi_{i=1}^{n+1} G_{\alpha_i}
 \rightarrow [X_1^{\alpha_1}:X_2^{\alpha_2}:\cdots: X_{n+1}^{\alpha_{n+1}}]
\in\Pn \alpha
$$

Let $d$ be a natural number.
The polynomial (resp. the polynomial form) $\omega$ in $\C^{n+1}$
is weighted homogeneous of degree $d$ if
$$
\lambda^*(\omega)=\lambda^d \omega, \ \lambda\in\C^*
$$
For a polynomial $g$ this means that
$$
g(\lambda^{\alpha_1}X_1,\lambda^{\alpha_2}X_2,\ldots, \lambda^{
\alpha_{n+1}}X_{n+1})=
\lambda^d g(X_1,X_2,\ldots,X_{n+1}),\ \forall \lambda\in\C^*
$$

Let $g$ be an irreducible polynomial of (weighted) degree $d$.
The set $g=0$ induces a hypersurface $D$ in $\Pn \alpha, \ \alpha=
(\alpha_1,\alpha_2,\ldots,\alpha_{n+1})$. If  $g$ has an isolated singularity 
at $0\in\C^{n+1}$ then Steenbrink has proved that  $D$ is a 
$V$-manifold/quasi-smooth variety. A $V$-manifold 
may be  singular  but it has many common features with smooth varieties (see 
\cite{st77,do82}).

For a polynomial form $\omega$ of degree $dk,\ k\in\N$
in $\C^{n+1}$ we have $\lambda^*\frac{\omega}{g^k}=\frac{\omega}{g^k}$ for all
$\lambda\in\C^*$. Therefore, $\frac{\omega}{g^k}$ induce a meromorphic form
on $\Pn \alpha$ with poles of order $k$ along $D$. If there is no confusion
we denote it again by $\frac{\omega}{g^k}$. The polynomial form
\begin{equation}
\label{23may04}
\eta_\alpha=\sum_{i=1}^{n+1}(-1)^{i-1}\alpha_iX_i\widehat {dX_i}
\end{equation}
where $\widehat{dX_i}=dX_1\wedge\cdots\wedge
 dX_{i-1}\wedge dX_{i+1}\wedge\cdots\wedge dX_{n+1}$, is of degree
$\sum_{i=1}^{n+1} \alpha_i$.

Let $\Pn {(1,\alpha)}=
\{[X_0:X_1:\cdots:X_{n+1}]\mid (X_0,X_1,\cdots,X_{n+1})\in
\C^{n+2}\}$ be the projective space of weight $(1,\alpha),\
\alpha=(\alpha_1,\ldots,\alpha_{n+1})$.
One can consider $\Pn {(1,\alpha)}$  as a compactification of
$\C^{n+1}=\{(x_1,x_2,\ldots,x_{n+1})\}$ by putting
\begin{equation}
\label{bala}
x_i=\frac{X_i}{X_0^{\alpha_i}},\ i=1,2,\cdots,n+1
\end{equation}
The projective space at infinity
$\Pn {\alpha} _{\infty}=\Pn {(1,\alpha)}-\C^{n+1}$ is of weight
$\alpha:=(\alpha_1,\alpha_2,\ldots,\alpha_{n+1})$.

Let $f$ be the strongly tame polynomial of (weighted) degree $d$ in the introduction
and $g$ be its last
quasi-homogeneous part. Now we can look at $f$ as
a rational function on $\Pn {(1,\alpha)}$ and the fibration $f=t$ as a pencil
in $\Pn {(1,\alpha)}$ with the axis $\{g=0\}\subset \Pn \alpha_{\infty}$. Note
that $\Pn \alpha _{\infty}$ itself is a fiber of this pencil.
This implies that the closure ${\overline L_c}$ of
 $L_c:=f^{-1}(c)$  in $\Pn {(1,\alpha)}$ 
intersects $\Pn \alpha _\infty$  transversally in the sense of $V$-manifolds.
In particular, 1. $f$ has connected fibers because $f$ at the infinity has 
connected fibers 2. $f$  has only isolated 
singularities because every algebraic
variety of dimension greater than zero in $\Pn {(1,\alpha)}$ intersects
$\Pn \alpha _\infty$. After making a blow-up along the axis
$\{g=0\}\subset \Pn \alpha_\infty$ and using Ehresmann's fibration theorem
one concludes that $f$ is $C^\infty$ fiber bundle over $\C\backslash C$.
\section{Hodge cycles}
\label{hodge}
Let $M$ be a smooth projective complex  manifold of
dimension $n$. The
cohomologies of $M$ with complex coefficients  carry the so called Hodge
decomposition
\begin{equation}
\label{11dez02}
H^m(M,\C)=H^{m,0}\oplus H^{m-1,1}\oplus\cdots\oplus H^{1,m-1}\oplus H^{0,m}
\end{equation}
Using de Rham cohomology
$$
H^m(M,\C)\cong H^{m}_{deR}(M):=
\frac{Z^m_d}{dA^{m-1}}
$$
 we have $H^{p,q}\cong\frac{Z^{p,q}_d}{dA^{p+q-1}\cap Z^{p,q}_d}$,
where  $A^{m}$ (resp. $Z^{m}_d$, $Z^{p,q}_d$)  is the set of $C^\infty$
differential $m$-forms (resp. closed  $m$-forms, closed $(p,q)$-forms)
on $M$ (with this notation one has  the canonical inclusions
$H^{p,q}\rightarrow H^m(M,\C)$ and one can  prove ~(\ref{11dez02}) using
harmonic forms, see M. Green's lectures \cite{gmv}, p. 14). The Hodge
filtration is defined
$$
G^p:=F^{p}H^m(M,\C)=H^{m,0}\oplus H^{m-1,1}\oplus\cdots\oplus H^{p,m-p}
$$
Let $m$ be an even natural number and $Z=\sum_{i=1}^s r_iZ_i$, where 
$Z_i,\ i=1,2,\ldots,s$ is a subvariety of $M$ of
complex dimension $\frac{m}{2}$ and $r_i\in \Z$. Using a
resolution map $\tilde Z_i\rightarrow M$, where $\tilde Z_i$ is a complex
manifold,  one can define an element
$\sum_{i=1}^s r_i[Z_i]\in H_m(M,\Z)$  which is called an algebraic cycle (see 
\cite{boha}).
 Since the restriction to $Z$ of a $(p,q)$-form with
$p+q=m$ and $p\not=\frac{m}{2}$ is identically zero, an algebraic cycle
$\delta$ has the following property:
$$
\int_\delta G^{\frac{m}{2}+1}=0
$$
 A cycle $\delta\in H_m(M,\Z)$  with the above property
is called a Hodge cycle. The assertion of the Hodge conjecture is that if
we consider the rational homologies then a Hodge cycle
$\delta\in H_n(M,\Q)$ is an algebraic
cycle, i.e. there exist subvarieties $Z_i\subset M$ of dimension
$\frac{m}{2}$ and rational numbers $r_i$  such that $\delta=
\sum r_i[Z_i]$.  The difficulty of this conjecture
lies in constructing varieties just with their homological information.

Now let $U$ be a quasi-projective smooth variety, $U\subset M$ its
compactification in the projective variety $M$ such that $N:=M-U$ is a
divisor with normal crossings (see \cite{de71} 3.2) and
$$
i:H_m(U,\Z)\rightarrow
 H_m(M,\Z)
$$
be the map induced by the inclusion $U\subset M$. For instance, 
$L_c=f^{-1}(c),c\in \C\backslash C$ of the previous section with $m=n$
is  an example of such a quasi-projective smooth
variety whose compactification divisor has only one irreducible 
component.  We are
 interested
to identify Hodge cycles in the image of $i$.
\begin{rem}\rm
\label{13.12.02}
Let $M$ be a hypersurface of even dimension
$n$ in the projective space $\Pn {n+1}$.
By the first Lefschetz theorem $H_m(M,\Z)\cong
H_m(\Pn{n+1},\Z)$, $m<n$ and so the only interesting Hodge cycles are
in $H_n(M,\Z)$. For a general  hypeplane section $N$ of $M$, 
the long exact sequence of the pair $(M,U)$, where $U$ is the
complement of $N$ in $M$, gives rise to the isomorphism
$H_{{\rm prim}}^n(M,\Q)\cong W_nH^n(U,\Q)$ induced by the 
inclusion $U\subset M$.
This implies that we capture all Hodge cycles in the primitive
cohomology (using Poincar\'e duality).    


\end{rem}
The mixed Hodge structure of
$H^m(U,\Q)$ consists of two filtrations
\begin{equation}
\label{21503}
 0=F^{m+1}\subset F^m\subset \cdots \subset F^1
\subset F^0=H^m(U,\C)
\end{equation}
$$
0=W_{m-1}\subset W_m\subset W_{m+1}\subset \cdots W_{m+a-1}\subset W_{m+a}
=H^m(U,\C), 1\leq a\leq m
$$
where  $W$  is defined over $\Q$, 
i.e. it is defined on $H^n(U, \Q)$ and we have 
complexified it. The number $a$ is the minimum of $m$ and the number of
irreducible components of $N$. Therefore, it is $1$ for $L_c$.
The first is the Hodge filtration and the second is the weight filtration. 
The Hodge filtration induces a filtration on 
$Gr_a^W:=W_a/W_{a-1}$ and we set 
\begin{equation}
\label{18.5.2004}
Gr_F^bGr^W_a:= F^{b}Gr_a^W/F^{b+1}Gr_a^W=
\frac{(F^{b}\cap W_a)+ W_{a-1}}{(F^{b+1}\cap W_a)
+ W_{a-1}}, a,b\in \Z 
\end{equation}
Let $r: H^m(M,\C)\rightarrow H^m(U,\C)$ be induced by inclusion. 
We have $W_m=r(H^m(M,\C))$ and $r$ is a weight zero morphism of mixed Hodge 
structures. It is strict and in  particular
$$
r(G^{\frac{m}{2}+1})=
 F^{\frac{m}{2}+1}\cap Im(r)=
F^{\frac{m}{2}+1}\cap W_m
$$
(see for instance \cite{kk} for definitions).
Now let us  be given  a cycle $\delta\in H_m(U,\Z)$ whose image in
$H_m(M,\Z)$ is Hodge.
The condition of being a Hodge cycle
translate
into a property of $\delta$ using the mixed Hodge structure of
$H^m(U,\C)$ as follows
$$
\int_{i(\delta)}G^{\frac{m}{2}+1}=\int_\delta r(G^{\frac{m}{2}+1})=
$$
\begin{equation}
\label{21.5.03}
\int_\delta F^{\frac{m}{2}+1}\cap W_{m}=0
\end{equation}
\begin{defi}
\label{20.5.2004}
A cycle $\delta\in H_m(U,\Z)$ is called Hodge
if (\ref{21.5.03}) holds, where $F^{\frac{m}{2}+1}$ (resp. $W_m$) is the
$(\frac{m}{2}+1)$-th (resp. $m$-th) piece of the Hodge filtration (resp. weight
filtration) of the mixed Hodge structure of $H^m(U,\Q)$.
\end{defi}
All the elements in the kernel of $i$ are Hodge cycles and we call
them trivial Hodge cycles. 
\section{Global Brieskorn modules}
\label{brieskorn}
In this section we introduce two Brieskorn modules $H'$ and $H''$ associated
to a polynomial $f$ and the notion of Gauss-Manin connection on them. In the
usual definition of Gauss-Manin connection for $n$-th cohomology of
the fibers of $f$,  if we take global sections and then compose it with
the vector field $\frac{\partial}{\partial t}$ in $\C$
then we obtain our version  of Gauss-Manin connection.

Let $f$ be the strongly tame polynomial in the introduction.
Multiplying by $f$ defines a linear operator on
\begin{equation}
\label{8.3.4}
V_f:=\frac{\C[x]}{<\frac{\partial f}{\partial x_i}\mid i=1,2,\ldots,n+1>}
\end{equation}
which is denoted by $A$. In the previous section we have seen that
$f$ has isolated
singularities and so  $V_f$ is a $\C$-vector space of finite dimension
$\mu$, where $\mu$ is the sum of local Milnor numbers of $f$,  and
eigenvalues of $A$ are exactly the critical values of $f$ (see for
instance \cite{mov}, Lemma 1.1).
Let  $S(t)\in\C[t]$ be the minimal polynomial of $A$, i.e. the polynomial
with the  minimum degree and with the leading coefficient $1$ such that
$S(A)\equiv 0$ as a function from $V_f$ to $V_f$
$$
S(f)=\sum_{i=1}^{n+1}p_i\frac{\partial f}{\partial x_i}, \ p_i\in\C[x]
$$
or equivalently
\begin{equation}
\label{alirezajan}
S(f)dx=df\wedge\eta_f, \ \eta_f=\sum_{i=1}^{n+1}(-1)^{i-1}p_i\widehat{dx_i}
\end{equation}
From now on we fix an $\eta_f$ with
the above property. To calculate $S(t)$ we may start with the characteristic polynomial   
$S(t)=det(A-tI)$,  where $I$ is the
$\mu\times \mu$ identity matrix  and  we have fixed a monomial 
basis of $V_f$ and have written $A$ as a matrix. This $S$ has the property
(\ref{alirezajan}) but it is in general useless from computational point of view 
(see \S \ref{examples}).
Note that if $f$ has rational coefficients then $S$ has 
also rational coefficients. The polynomial 
$S$ has only zeros at critical values  $C$ 
of $f$.
 
The global Brieskorn modules are
$$
{H''}=\frac{\Omega^{n+1}}{df\wedge d\Omega^{n-1}}, \
H'=\frac{df\wedge \Omega^n}{df\wedge d\Omega^{n-1}}
$$
They are $\C[t]$-modules. Multiplication by $t$ corresponds
to the usual multiplication of differential forms with $f$.
The Gauss-Manin connection
$$
\nabla: {H'}\rightarrow {H''},\ \nabla ([df\wedge \omega])=[d\omega]
$$
is a well-defined function  and satisfies the Leibniz rule
\begin{equation}
\label{shaban}
\nabla (p\omega)=p\nabla(\omega)+p'\omega,\ p\in\C[t],\ \omega\in H'
\end{equation}
where $p'$ is the derivation with respect to $t$.
Let $H_C'$ (resp. $H_C''$ and $\C[t]_C$) be the localization of
$H'$ (resp. $H''$ and $\C[t]$) on the multiplicative
subgroup of $\C[t]$ generated by $\{t-c,c\in C\}$.
An element of $H_C'$ is a fraction
$\omega/p,\ \omega\in H',\ p\in\C[t], \{p=0\}\subset C$.
Two such fractions
$\omega/p$ and $\tilde\omega/\tilde p$ are equal if
 $\tilde p\omega=p\tilde\omega$.
We have  $\frac{H''}{H'}=V_f$ and
so $S.{H''}\subset {H'}$. This means that the inclusion $H'\subset H''$
induces an equality $H_C'=H_C''$.
We denote by  $H_C$ the both side of the equality.
Let
$\tilde \Omega^i$ denote the set of rational differential $i$- forms
in  $\C^{n+1}$ with poles along the $L_c,\ c\in C$.
The canonical map $H_C\rightarrow
\frac{\tilde \Omega^{n+1}}{df\wedge d\tilde \Omega^{n}}$ is an isomorphism
of $\C[t]_C$-modules and this gives another interpretation of $H_C$.
 One extends $\nabla$ as a function
from $H_C$ to itself by
\begin{equation}
\label{12may04}
\nabla(\frac{[df\wedge\omega]}{p})=[d(\frac{\omega}{p(f)})]=
\frac{p[d\omega]- p'[df\wedge\omega]}{p^2}, \ p\in\C[t],\
[df\wedge\omega]\in H'
\end{equation}
This is a natural extension of $\nabla$ because
it satisfies
$$
\nabla (\frac{\omega}{p})=\frac{p\nabla\omega-p'\omega}{p^2},\ p\in\C[t],
\ \omega \in H_C
$$
\begin{lemm}
We have
$$
\nabla([Pdx])=\frac{[(Q_P-P.S'(f))dx]}{S}, \ P\in\C[x]
$$
where
$$
Q_P=\sum_{i=1}^{n+1}
(\frac{\partial P}{\partial x_i}p_i+P\frac{\partial p_i}{\partial x_i})
$$
\end{lemm}
\begin{proof}
\begin{eqnarray*}
\nabla([Pdx]) & =& \nabla (\frac{[df\wedge P\eta_f]}{S})=
[d(\frac{P\eta_f}{S(f)})]=
\frac{[S(f)d(P\eta_f)-S'(f)P df\wedge\eta_f]}{S^2} \\
 &=&
\frac{[d(P\eta_f)-S'(f)Pdx]}{S}=
\frac{[dP\wedge\eta_f+P.d\eta_f- P.S'(f)dx]}{S} \\
 &=&
\frac{[(Q_P-P.S'(f))dx]}{S}
\end{eqnarray*}
\end{proof}
It is better to have in mind that
the polynomial $Q_P$  is defined by the relation $d(P\eta_f)=Q_Pdx$.
In the next section we will use the iterations of Gauss-Manin connection,
$\nabla^k=\nabla\circ\nabla\circ\cdots\circ\nabla, k$ times. 
To be able to calculate them  we need the following
operators
$$
\nabla_k:H''\rightarrow H'',\
 k=0,1,2,\ldots
$$
$$
 \nabla_k(\omega)=
\nabla(\frac{\omega}{S(t)^k})S(t)
^{k+1}=
S(t)\nabla(\omega)-k.S'(t)\omega
$$
For $\omega=Pdx$ we obtain the formula
$$
\nabla_k(Pdx)=(Q_P-(k+1)S'(t)P)dx
$$
We show by induction on $k$ that
\begin{equation}
\nabla^k=
\frac{\nabla_{k-1}\circ\nabla_{k-2}\circ\cdots\circ\nabla_0}{S(t)^{k}}
\end{equation}
The case $k=1$ is trivial. If the equality is true for $k$ then
\begin{eqnarray*}
\nabla^{k+1} &= & \nabla\circ\nabla^{k}=
\nabla(\frac{\nabla_{k-1}\circ\nabla_{k-2}\circ\cdots\circ\nabla_0}{S(t)^{k}})
\\
 &= &
\frac{\nabla_0\circ\nabla_{k-1}\circ\nabla_{k-2}\circ\cdots\circ\nabla_0-
kS'(t)\nabla_{k-1}\circ\nabla_{k-2}\circ\cdots\circ\nabla_0}{S(t)^{k+1}} \\
 &=&
\frac{\nabla_{k}\circ\nabla_{k-1}\circ\cdots\circ\nabla_0}{S(t)^{k+1}}
\end{eqnarray*}
The Brieskorn module $H'=
\frac{\Omega^n}{df\wedge\Omega^{n-1}+d\Omega^{n-1}}$ defined
in the introduction is isomorphic to the one in this section
by the map $[\omega]\rightarrow [df\wedge\omega]$. The inverse of the canonical
isomorphism $H_C'\rightarrow H_C''$ is denoted by
$\omega\in H_C''\rightarrow \frac{\omega}{df}\in H_C'$.
The Gauss-Manin connection with this notation
can be written in the form
$$
\nabla: H'\rightarrow H_C',\ \nabla(\omega)=\frac{d\omega}{df}:=
\frac{\omega_1}{S(t)}
$$
where $S(f)d\omega=df\wedge\omega_1$. In the literature 
one also calls $\frac{d\omega}{df}$ the Gelfand-Leray form of $d\omega$. 
Looking in this way it
turns out that
\begin{equation}
\label{mohebbi}
df\wedge \nabla\omega=\nabla(df\wedge\omega),\ \forall \omega\in H'
\end{equation}
Let $U$ be an small open set in $\C\backslash C$, 
$\delta_t\in H_n(L_t,\Z),\ t\in U$ be a continuous family of cycles 
 and $\omega\in H'$.
The main property of the Gauss-Manin connection is
\begin{equation}
\label{20.5.04}
\frac{\partial}{\partial t}\int_{\delta_t}\omega=\int_{\delta_t}
\nabla\omega
\end{equation}
Recall the notations introduced for a quasi-homogeneous polynomial $f=g$ 
in Introduction. For this $f$
$S(t)=t$ and $\eta_f$ is $\eta$ in (\ref{21may2004}).
This means that $fdx=df\wedge \eta$.
Since this equality is linear in $f$ it is enough to check it for
monomials $x^\alpha,\alpha.w=1$.
$$
dx^\alpha\wedge\eta= (\sum_{i=1}^{n+1}\alpha_i\frac{x^\alpha}{x_i}dx_i)
\wedge(\sum_{i=1}^{n+1}(-1)^{i-1}
w_ix_i\widehat {dx_i})=(\alpha.w) x^\alpha dx=x^\alpha dx
$$
We have also  $d\eta= (w.1) dx$.
$$
d\eta_\beta=dx^\beta\wedge\eta+ x^\beta d\eta=(\beta.w) x^\beta dx+
(w.1)x^\beta dx=A_\beta x^\beta dx=\frac{A_\beta }{f}df\wedge(x^{\beta}\eta)
$$
which implies that $\nabla\eta_\beta=\frac{A_\beta}{t}\eta_\beta$
( In the same way on can check that
$\nabla(x^\beta dx)=\frac{(A_\beta-1)}{t}x^\beta dx$).
This implies that
$\frac{\partial}{\partial t}\int_{\delta_t}\eta_\beta=\frac{A_\beta}{t}
\int_{\delta_t}\eta_\beta$.
Therefore there exists a constant number $C$ depending only on $\eta_\beta$ 
and $\delta_t$ such that $\int_{\delta_t}\eta_\beta=Ct^{A_\beta}$. 
One can take a branch of $t^{A_\beta}$ so that $C=\int_{\delta_1}\eta_\beta$.
\section{Residue map on the Brieskorn module}
\label{residue}
\def\res{{\mathrm Res}}
Let us be given a closed submanifold $N$ of
pure real codimension $c$ in  a manifold $M$.
The  Leray (or Thom-Gysin) isomorphism is
$$
\tau :H_{k-c}(N,\Z)\tilde{\rightarrow}H_k(M,M- N,\Z)
$$
holding for any $k$, with the convention that $H_s(N)=0$ for
$s<0$. Roughly speaking, given $x\in H_{k-c}(N)$, its image by this
isomorphism is obtained by thickening a cycle representing $x$,
  each point of it growing into a closed $c$-disk transverse to $N$ in $M$
  (see  for instance \cite{che} p. 537).
Let $N$ be a connected codimension one submanifold of the complex manifold
  $M$ of dimension $n$.
Writing the long exact sequence of the pair $(M,M-N)$ and using
$\tau$ we obtain:
\begin{equation}
\label{longler}
\cdots\rightarrow H_{n+1}(M,\Z)\rightarrow H_{n-1}(N,\Z)
\stackrel{\sigma}{\rightarrow}H_n(M-N,\Z)
\stackrel{i}{\rightarrow} H_n(M,\Z)\rightarrow \cdots
\end{equation}
where $\sigma$ is the composition of the boundary operator with
$\tau$ and $i$ is induced by inclusion. Let
$\omega\in H^n(M-N,\C):=H_n(M-N,\Z)^*\otimes\C$, where $H_n(M-N,\Z)^*$ is the
dual of $H_n(M-N,\Z)$ . The composition
$\omega\circ \sigma: H_{n-1}(N,\Z)\rightarrow \C$ defines a linear map
and its complexification is an element in $H^{n-1}(N,\C)$. It is denoted
by $\res _{N}(\omega)$ and called the residue of $\omega$ in $N$.
We consider the case in which $\omega$ in the
$n$-th de Rham cohomology of $M-N$ is represented by a meromorphic
$C^\infty$ differential form $\omega'$ in $M$ with poles of order at most one
along $N$.
Let $f_\alpha=0$ be the defining equation of $N$ in a neighborhood $U_\alpha$
of a  point $p\in N$ in $M$ and write $\omega'=\omega_\alpha
\wedge\frac{df}{f}$. For two such neighborhoods $U_\alpha$ and $U_\beta$
with non empty intersection we have
$\omega_\alpha=\omega_\beta$ restricted to $N$.
Therefore we get a $(n-1)$-form on $N$
which in the de Rham cohomology of $N$ represents $\res_N{\omega}$
(see \cite{gr69} for details). This is called Poincar\'e residue.
The residue map $\res_N$ is a morphism
of weight $-2$ of mixed Hodge structures, i.e. 
$$
\res_N(W_pH^n(M-N,\C))\subset\res_N(W_{p-2}H^{n-1}(N,\C)),\ p\in\Z
$$
$$
\res_N(F^qH^n(M-N,\C))\subset\res_N(F^{q-1}H^{n-1}(N,\C)),\ q\in\Z 
$$
We fix a regular value $c\in\C\backslash C$.
To each $\omega\in H_C''$ we can associate the
residue of $\frac{\omega}{(t-c)^k}$ in $L_c$ which is going to be
an element of $H^n(L_c,\C)$. This map is well-defined because
$$
\frac{df\wedge d\omega}{p(f)(f-c)^k}=d(\frac{df\wedge\omega}{p(f)(f-c)^k}),\
p\in\C[t],\ \omega\in \Omega^{n-1}
$$

\begin{lemm}
\label{hayhay}
For $\omega\in H''$ and $k=2,3,\ldots$ the forms $\frac{\omega}{(t-c)^k}$ and
$\frac{\nabla\omega}{(k-1)(t-c)^{k-1}}$ have the same residue in $L_c$.
In particular the residue of $\frac{\omega}{(t-c)^k}$  in $L_c$ is
the restriction of $\frac{\nabla^{k-1}\omega}{(k-1)!df}$ to $L_c$ and
the residue of $\frac{df\wedge \omega}{(t-c)^k},\ \omega\in H'$
in $L_c$ is the restriction of $\frac{\nabla^{k-1}\omega}{(k-1)!}$ to $L_c$.
\end{lemm}
\begin{proof}
This Lemma is well-known in the theory of Gauss-Manin systems 
(see for instance \cite{st85}). We give
an alternative proof in the context of this article.
We have
$$
\nabla(\frac{\omega}{(t-c)^{k-1}})=\frac{\nabla(\omega)}{(t-c)^{k-1}}-(k-1)
\frac{\omega}{(t-c)^k}
$$
According to (\ref{12may04}),
the left hand side corresponds to an exact form and so it
does not produce a residue in $L_c$. This proves that first
part. To obtain the second part we repeat $k-1$ times the result
of the first part on $\frac{\omega}{(t-c)^k}$ and we obtain
$\frac{\nabla^{k-1}\omega}{(k-1)!(t-c)}$.
Now we take the Poincar\'e residue and
obtain the second statement. The third statement is a consequence of
the second and the identity (\ref{mohebbi}).
\end{proof}
Note that the residue of $\frac{\omega}{t-c},\ \omega\in H''$ in $L_c$ 
coincides with the restriction of $\frac{\omega}{df}\in H_C'$ to $L_c$.
\section{Griffiths-Steenbrink Theorem}
\label{grst}
This section is dedicated to a classic theorem of Griffiths 
in \cite{gr69}. Its generalization for quasi-homogeneous spaces is due
to Steenbrink in \cite{st77}. In both cases there is not given 
an explicit basis of the Hodge structure of the complement of a smooth
hypersurface. This is the main reason to put Theorem \ref{1.2.04}  
in this article. Recall the notations of
\S \ref{wps}.
\begin{lemm}
\label{2004}
For a monomial $x^\beta$ with $A_\beta=k\in \N$,
the meromorphic form $\frac{x^\beta dx}{(f-t)^k}$
has a pole of order one at infinity and its Poincar\'e residue at infinity is
$\frac{X^\beta\eta_\alpha}{g^k}$.
 \end{lemm}
\begin{proof}
Let us write the above form in the homogeneous coordinates (\ref{bala}).
We use $d(\frac{X_i}{X_0^{\alpha_i}})=X_0^{-\alpha_i}dX_i-\alpha_iX_iX_0^{-\alpha_i-1}dX_0 $ and
\begin{eqnarray*}
\frac{x^\beta dx}{(f-t)^k} &= &
\frac
{(\frac{X_1}{X_0^{\alpha_1}})^{\beta_1}\cdots(\frac{X_{n+1}}{X_0^{\alpha_{n+1}}})
^{\beta_{n+1}}d(\frac{X_1}{X_0^{\alpha_1}})\wedge\cdots
\wedge d(\frac{X_{n+1}}{X_0^{\alpha_{n+1}}})}
{(f(\frac{X_1}{X_0^{\alpha_1}},\cdots,\frac{X_{n+1}}{X_0^{\alpha_{n+1}}})-t)^k}
\\
 &=&
\frac
{X^\beta\eta_{(1,\alpha)}}
{X_0^{(\sum_{i=1}^{n+1} \beta_i\alpha_i)+(\sum_{i=1}^{n+1}\alpha_i)+1-kd}
(X_0\tilde F-g(X_1,X_2,\cdots,X_{n+1}))^k} \\
&= &
\frac{X^\beta\eta_{(1,\alpha)}}{X_0(X_0\tilde F-g(X_1,X_2,\cdots,X_{n+1}))^k}
\\
&=&
\frac{dX_0}{X_0}\wedge \frac{X^\beta\eta_\alpha}{(X_0\tilde F-g)^k}
\end{eqnarray*}
The last equality is up to forms without pole at $X_0=0$.
The restriction of $\frac{X^\beta\eta_\alpha}{(X_0\tilde F-g)^k}$
to $X_0=0$ gives
us the desired form.
\end{proof}

\begin{theo}(Griffiths-Steenbrink)
\label{1.2.04}
Let $g(X_1,X_2,\cdots,X_{n+1})$ be a weighted
homogeneous polynomial of degree $d$, weight
$\alpha=(\alpha_1,\alpha_2,\ldots,\alpha_{n+1})$ and with an isolated
singularity at $0\in\C^{n+1}$(and so
$D=\{g=0\}$ is a $V$-manifold).  We have
$$
H^n(\Pn \alpha-D,\C)\cong\frac{H^0(\Pn\alpha, \Omega^n(*D))}
{dH^0(\Pn \alpha,\Omega^{n-1}(*D))}
$$
and under the above isomorphism
\begin{equation}
\label{bulu}
Gr_F^{n+1-k}Gr^W_{n+1}H^{n}(\Pn \alpha-D,\C):=F^{n-k+1}/F^{n-k+2}\cong
\end{equation}
$$
\frac{H^0(\Pn\alpha, \Omega^n(kD))}
{dH^0(\Pn \alpha,\Omega^{n-1}((k-1)D))+H^0(\Pn \alpha,\Omega^{n}((k-1)D))}
$$
where $0=F^{n+1}\subset F^{n}\subset\cdots\subset F^1\subset F^0=
H^n(\Pn \alpha-D,\C)$ is the Hodge filtration of $H^n(\Pn \alpha-D,\C)$.
Let $\{X^\beta\mid\beta\in I\}$ be a basis of monomials for the
Milnor vector space
$$
\C[X_1,X_2,\cdots,X_{n+1}]/
<\frac{\partial g}{\partial X_i}\mid i=1,2,\ldots, n+1>
$$
A basis of (\ref{bulu}) is given by
\begin{equation}
\label{27.1.04}
\frac{X^\beta\eta_\alpha}{g^k},\ \beta\in I,\ A_\beta=k
\end{equation}
where $\eta_\alpha$ is given by (\ref{23may04}).
\end{theo}
Recall that if $D$ is normal crossing divisor in a projective
variety $M$ then $H^m(M-D,\C)\cong {\mathbb H}^m(M,\Omega^{\bullet}(\log D)),
m\geq 1$ and the $i$-th piece of the Hodge filtration of $H^m(M-D,\C)$ under
this isomorphism is given by
${\mathbb H}^m(M,\Omega^{\bullet\geq i}(\log D))$ (see \cite{de71}). By
definition we have $F^0/F^1\cong H^{n}(M,\O_M)$ and so in the situation
of the above theorem $F^0=F^1$. Note that in the above theorem the
 residue map $r:H^n(\Pn \alpha-D,\C)\rightarrow H^{n-1}(D,\C)_0$ is an
isomorphism of Hodge structures of weight $-2$, i.e. it maps
the $k$-th piece of the Hodge filtration of  $H^n(\Pn \alpha-D,\C) $ to
 the $(k-1)$-th piece of the Hodge filtration of $ H^{n-1}(D,\C)_0$. Here
the sub index $0$ means the primitive cohomology.
\begin{proof}

The first part of this theorem for usual projective spaces is
due to Griffiths \cite{gr69}. The generalization for quasi-homogeneous spaces
is due to Steenbrink \cite{st77}. The essential ingredient in the proof
is Bott's vanishing theorem for quasi-homogeneous spaces:
Let $L\in H^1(\Pn \alpha ,{\cal O}^*)$ be a line bundle on $\Pn \alpha$ with $c(\pi^{*}L)=k$, 
where $\Pn n\rightarrow \Pn \alpha$ is the canonical map. Then
$H^p(\Pn \alpha,\Omega^q_{\Pn \alpha}\otimes L)=0$ except possibly in the case $p=q$ and $k=0$, or
$p=0$ and $k>q$, or $p=n$ and $k<q-b$.

The proof of the second part which gives an explicit basis of
Hodge filtration is as follows: We consider $\Pn \alpha$ as the
projective space at infinity in $\Pn {(1,\alpha)}$.
According to Lemma \ref{2004} for $f=g$ the residue of the form
$\frac{x^\beta dx}{(g-1)^k}$ at $g=1$ is (\ref{27.1.04}). Now we use
Lemma 5 of \cite{st77}. This lemma says that the residue of
$\frac{x^\beta dx}{(g-1)^k}$ at infinity form a basis of (\ref{bulu}).
\end{proof}

\section{An explicit basis of $H_{\tilde C}$}
\label{basis}
In this section we prove Theorem \ref{main}. 
In the following by homogeneous
we mean weighted homogeneous with respect to fixed weights
$\alpha=(\alpha_1,\alpha_2,\ldots,\alpha_{n+1})$.

Let $f$ be the polynomial in the introduction, $f=f_0+f_1+f_2+\cdots+
f_{d-1}+f_d$ be its homogeneous decomposition in the graded ring
$\C[x_1,x_2,\ldots,x_{n+1}],\ \deg(x_i)=\alpha_i$ and $g:=f_d$ the last 
homogeneous part of $f$. 
Let also  $F=f_0x_0^d+f_1x_0^{d-1}+\cdots+f_{d-1}x_0+g$ be the
homogenization of $f$.
\begin{lemm}
\label{2.10.04}
The set $x^I$ generates freely the $\C[x_0]$-module
$$
V:=\C[x_0,x]/<\frac{\partial F}{\partial  x_i}\mid i=1,2,\ldots,n+1>
$$
\end{lemm}
\begin{proof}
First we prove that $x^I$ generates  $V$ as a $\C[x_0]$-module.
We write the expansion of $P(x_0,x)=\sum_{i=0}^kx_0^iP_i(x)\in\C[x_0,x]$
in $x_0$ and so it is enough to prove that every element $P\in \C[x]$ can be
written in the form
\begin{equation}
\label{3.3.04}
P=\sum_{\beta\in I}C_\beta x^\beta+\sum_{i=1}^{n+1}
Q_i\frac{\partial F}{\partial x_i},\ C_\beta \in\C[x_0],\  Q_i\in\C[x_0,x]
\end{equation}
Since $x^I$ is a basis of $V_g$, we can write
\begin{equation}
\label{divane}
 P=\sum_{\beta\in I}c_\beta x^\beta+\sum_{i=1}^{n+1}
q_i\frac{\partial g}{\partial x_i},\ c_\beta \in\C,\  q_i\in\C[x]
\end{equation}
We can choose $q_i$'s so that
\begin{equation}
\label{divan}
\deg(q_i)+\deg(\frac{\partial g}{\partial x_i})\leq \deg (P)
\end{equation}
If this is not the case then we write the non-trivial
homogeneous equation of highest degree obtained from (\ref{divane}). Note
that $\frac{\partial g}{\partial x_i}$ is homogeneous.
If some terms  of $P$ occur in this new equation then
we have already (\ref{divan}). If not we subtract this new equation from
(\ref{divane}). We repeat this until getting the first case and so the
desired inequality.
Now we have
$\frac{\partial g}{\partial x_i}=
\frac{\partial F}{\partial x_i}-x_0\sum_{j=0}^{d-1}
\frac{\partial f_j}{\partial x_i}x_0^{d-j-1}$ and so
\begin{equation}
\label{divane1}
 P=\sum_{\beta\in I}c_\beta x^\beta+\sum_{i=1}^{n+1}
q_i\frac{\partial F}{\partial x_i}-x_0(\sum_{i=1}^{n+1}
\sum_{j=0}^{d-1}
q_i\frac{\partial f_j}{\partial x_i}x_0^{d-j-1})
\end{equation}
From (\ref{divan}) we have $(*): \deg(q_i\frac{\partial f_j}{\partial x_i})
\leq \deg(P)-1$. We write again $q_i\frac{\partial f_j}{\partial x_i}$ in the
form (\ref{divane}) and substitute it in (\ref{divane1}). By degree conditions
this process stops and at the end we get the equation (\ref{3.3.04}).

Now let us prove that $x^I$ generate the $\C[x_0]$-module freely.
For every $x_0=a$ fixed, let $V_{a}$ be the specialization of $V$ at
$x_0=a$. All $V_a$'s are vector spaces of the same dimension and according
to the above argument $x^I$ generates all $V_a$'s.
For $V_0$ it is even a $\C$-basis and so $x^I$ is a basis of all $V_a$'s.
If  the elements of $x^I$ are not $\C[x_0]$-independent then
we have $C.\sum_{\beta\in I}{C_\beta}x^\beta=0$ in $V$ for
some $C,C_\beta\in \C[x_0]$
and $C_\beta$'s do not have common zeros.
 We take an $a$ which is not a zero of $C$.
We have $\sum_{\beta\in I}{C_\beta(a)}x^\beta=0$ in $V_a$ which
is a contradiction.
\end{proof}

\begin{prop}
\label{28.6.04}
For every strongly tame polynomial $f\in\C[x]$ the forms 
$\omega_\beta:=x^\beta dx,\beta\in I$ (resp. $\eta_\beta:=x^\beta\eta,\ 
\beta\in I$) form a $\C[t]$- basis of the Brieskorn module $H''$  (resp. $H'$) of $f$.  
\end{prop}
\begin{proof}
We first prove the statement for $H''$. The statement for $f=g$ is well-known
(see for instance \cite{arn}). Recall the definition of the degree of a 
form in \S \ref{wps}. 
We write an element $\omega\in \Omega^{n+1}, \deg(\omega)
=m$ 
in the form 
$$
\omega=\sum_{\beta\in I} p_\beta(g)\omega_\beta+dg\wedge 
d\psi,\ p_\beta\in\C[t],\ \psi\in \Omega^{n-1},\ 
\deg(p_\beta(g) \omega_\beta)\leq m,\ \deg(d\psi)\leq m-d  
$$
This is possible because $g$ is homogeneous.
Now, we write the above equality in the form
$$
\omega=\sum_{\beta\in I} p_\beta(f)\omega_\beta+df\wedge d\psi+\omega',
\hbox{ with } 
\omega'=\sum_{\beta\in I} (p_\beta(g)-p_\beta(f))\omega_\beta+
d(g-f)\wedge d\psi
$$
The degree of $\omega'$ is strictly less than $m$ and so we repeat
what we have done at the beginning and finally we write $\omega$ as a 
$\C[t]$-linear combination of $\omega_\beta$'s. The forms $\omega_\beta, 
\beta\in I$ are linearly independent because $\# I=\mu$ and $\mu$ is 
the dimension
of $H^n(L_c,\C)$ for a regular $c\in\C-C$. The proof for $H'$ is similar and 
uses the fact that for $\eta\in\Omega^{n}$ one can write
\begin{equation}
\label{29.6}
\eta=\sum_{\beta\in I}p_\beta(g)\eta_\beta+dg\wedge\psi_1+d\psi_2 
\end{equation}
and each piece in the right hand side of the above equality has degree less 
than $\deg(\eta)$.

\end{proof}
The above proposition gives us an algorithm to write every element of $H'$ 
(resp. $H''$) 
as a $\C[t]$-linear sum of $\eta_\beta$'s (resp. $\omega_\beta$'s).
We must find such an algorithm first for the case $f=g$, which is not hard to
do (see \cite{mov04}). Note that if $\eta\in \Omega^{n}$ is written in the
form (\ref{29.6}) then
$$
d\eta=\sum_{\beta\in I}(p_\beta(g)A_\beta+p'_\beta(g)g)\omega_\beta-dg\wedge d\psi_1
$$

We specialize the module $V$ at $x_0=1$ and use
Lemma \ref{2.10.04} and obtain the following fact: $x^I$ form a basis for
the Milnor vector space $V_f$ of $f$. Let $F_t=F-t.x_0^d$.

\begin{lemm}
\label{10:35} Let $b\in\C\backslash C$. There is a map $\beta\in
I\rightarrow d_\beta\in \N\cup\{0\}$ such that the
$\C$-vector space $
\tilde V:=\C[x_0,x]/
<\frac{\partial F_b}{\partial  x_i}\mid i=0,1,\ldots,n+1>
$ is freely generated by
\begin{equation}
\label{14may04}
 \{ x_0^{\beta_0}x^\beta, 0\leq \beta_0\leq
d_\beta-1, \beta\in I\}
\end{equation}
In particular, the $\C(t)$-vector space
$$
 V':=\C(t)[x_0,x]/
<\frac{\partial F_t}{\partial  x_i}\mid i=0,1,\ldots,n+1>
$$
is freely generate by (\ref{14may04}).  
\end{lemm}
\begin{proof}
We consider the class $Cl$ of all sets of the form (\ref{14may04})
whose elements are linearly independent in $\tilde V$. For instance
the one element set $\{x^{\beta'}\}$ is in this class. 
In this example $d_{\beta'}=1$ and $d_\beta=0, \forall \beta\in I-
\{\beta'\}$.  Since $\tilde V$ is a finite dimensional $\C$-vector
space, $Cl$ has only a finite number of elements and so we can
take a maximal  element $A$ of $Cl$, i.e. there is no
element of $Cl$ containing $A$. We prove that $A$ generates $\tilde V$ 
and so it is the desired
set. Take a $\beta\in I$. We claim that $x^\beta x_0^k, \
k>d_\beta-1$, can be written as a linear combination of the elements
of $A$. The claim is proved by induction on $k$. For
$k=d_\beta$, it is true because $A$ is maximal and $A\cup
\{x^\beta x_0^{d_\beta}\}$ is of the form (\ref{14may04}).
Suppose that the claim is true for $k$. We write $x^\beta x_0^k$
as linear combination of elements of $A$ and multiply it by $x_0$.
Now we use the hypothesis of the induction for $k=d_\beta$ for
the elements in the new summand which are not in $A$ and get a
linear combination of the elements of $A$.
Since $\tilde V=V/<\frac{\partial F_b}{\partial x_0}>$ and $V$ is a
$C[x_0]$-module generate by $x^\beta$'s, we conclude that $x^\beta
x_0^k, k\in \N\cup\{0\},\ \beta\in I$ generate $\tilde V$ and so $A$
generates $\tilde V$.

If there is a $\C(t)$-linear relation between the elements of 
(\ref{14may04}) then we multiply it by a suitable element of $\C(t)$
and obtain a $\C[t]$-linear relation such that putting $t=b$ gives us
a nontrivial relation in $\tilde V$. This proves the second part.
\end{proof}
\begin{rem}\rm
\label{21may}
Lemma \ref{10:35} implies that for all $c\in \C$, except a
finite number of them which includes $C$ and does not include $b$, the set (\ref{14may04})
is a basis of the specialization of $V'$ at $t=c$. The set $\tilde C$ of
such exceptional values may be greater than $C$. To avoid such a
problem we may try to prove the following fact which seems to be
true: \\
(*) Let $\C[t]_C$ be the
localization of $\C[t]$ on its multiplicative subgroup generated
by $t-c,c\in C$. There is a function $\beta\in I\rightarrow
d_\beta\in \N\cup\{0\}$ such that the $\C[t]_C$-module
$V'':=
\C[t]_C[x_0,x]/<\frac{\partial F_t}{\partial  x_i}\mid i=0,1,\ldots,n+1>$ 
is freely generated by $\{
x_0^{\beta_0}x^\beta, 0\leq \beta_0\leq d_\beta-1, \beta\in I\}$.
\\
In \cite{mov04} I have used another algorithm 
(different from the one in the proof of Lemma \ref{10:35}). The advantage of 
this algorithm is that it also determines whether the obtained basis of $V'$ 
is a $\C[t]_C$ basis of $V''$ or not. A similar algorithm shows that one 
can take $d_\beta=d-1$ for a generic $b$.
\end{rem}
Let $f$ be a quasi-homogeneous polynomial of degree $d$. In this case 
$F_t=f-t.x_0^d$ and $V'$ is generated by $\{x^\beta x_0^{\beta_0},\beta\in I, 
0\leq \beta_0\leq d-2\}$. 
In fact it generates the $\C[t]_C$-module $V''$ freely and
so $\tilde C=C=\{0\}$.
The dimension of $V'$ is $(d-1)\mu$ and so
the Milnor number of $f-t.x_0^d$ is $(d-1)\mu$. Since the Milnor number
is topologically invariant, we conclude that for an arbitrary strongly tame function
$f$ the dimension of $V'$ is $(d-1)\mu$ and so
$$
\sum_{\beta\in I}d_\beta=(d-1)\mu
$$
Moreover,  we have $0<A_{(\beta,\beta_0)}=A_\beta+\frac{\beta_0+1}{d}<n+2$ for all $\beta\in I$ and 
$0 \leq \beta_0\leq d_\beta-1$ and so
$$
d_\beta<d(n+2-A_\beta)
$$  
{\it Proof of Theorem \ref{main}:}
Since the dimensions of the pieces of the mixed Hodge structure of a smooth fiber $L_c$ does not
depend on the analytic structure, the equality $Gr^W_mH'=0, \ m\not =n,n+1$ follows form
Steenbrink's theorem for the quasi-homogeneous polynomials. 
  
We use Lemma \ref{10:35} and we obtain a basis
$\{ x_0^{\beta_0}x^\beta, 0\leq \beta_0\leq d_\beta-1, \beta\in I\}$ of
the $\C[t]_{\tilde C}$-module $V''$. Recall the notations introduced in 
\S \ref{hodge} and $\tilde C$ in Remark \ref{21may}.
Theorem \ref{1.2.04} and Lemma \ref{2004} imply that the residue of the
forms
\begin{equation}
\label{darm}
\frac{x^\beta dx}{(f-c)^k},\ A_\beta=k
\end{equation}
in $L_c$ form a basis
of $Gr^{n+1-k}_FGr^W_{n+1}H^n(L_c,\C)$
(Residue map is morphism of weight -2 of mixed Hodge structures).
By Theorem  \ref{1.2.04} for $\Pn {(1,\alpha)}$ and the hypersurface
$X_c: F-cx_0^d=0, \ c\in \C\backslash \tilde C$ and Lemma \ref{10:35}
$$
\frac{x_0^{\beta_0}x^\beta\eta_{(1,\alpha)}}{(F-cx_0^d)^k}, \
A_\beta+\frac{\beta_0+1}{d}=k,\ 0 \leq \beta_0\leq d_\beta-1
$$
form a basis for  $Gr^{n+2-k}_FGr^W_{n+2} H^{n+1}(\Pn {(1,\alpha)}-X_c,\C)$.  
In the affine coordinate $\C^{n+1}\subset \Pn {(1,\alpha)}$, these
forms are
\begin{equation}
\label{darmstadt}
\frac{x^\beta dx}{(f-c)^k},\ 
A_\beta+\frac{1}{d}\leq k
\leq A_\beta+\frac{d_\beta}{d},\ k\in \N.
\end{equation}
So the residues of the above forms at $L_c$ form a basis  of
$Gr^{n+1-k}_F Gr^W_n H^{n}(L_c,\C)$.  We apply Lemma \ref{hayhay} to the
meromorphic forms (\ref{darm}) and (\ref{darmstadt}) and obtain the
fact that the forms (\ref{naucia}) (resp. the forms (\ref{vomiting}))
restricted to  the fiber $L_c,\ c\in\C\backslash \tilde C  $ form a basis of
$Gr^{n+1-k}_FGr_{n+1}^W H^n(L_c,\C )$ 
(resp. $Gr^{n+1-k}_FGr_{n}^W H^n(L_c,\C)$).
Note that $x^\beta dx=d(\frac{\eta_\beta}{A_\beta})$.

\qed
\section{Examples and applications}
\label{examples}
In this section we give some examples of the polynomial $f$ and 
discuss the result of the paper on them. The examples which 
we discuss are of the form 
$f(x_1,x_2,\ldots,x_{n+1})=\sum_{i=1}^{n+1}f_i(x_i)$,  where 
$f_i$ is a polynomial of degree $m_i, \  m_i\geq 2$ 
in one variable $x_i$ and with leading
coefficient one. 
Let $d$  be the least common multiple of  $m_i's$. We consider
$f$ in the weighted ring $\C[x], deg(x_i)=\frac{d}{m_i}, \ i=1,2\ldots, 
n+1$.
Then $deg(f)=d$ and the last homogeneous part of $f$ is
$g=x_1^{m_1}+x_2^{m_2}+\cdots+x_{n+1}^{m_{n+1}}$.
The vector space $V=\C[x]/ \hbox{{\rm Jacob}} (f)$ has the following basis of 
monomials 
$$
x^\beta,\ \beta\in I:=\{ \beta\in \Z^{n+1}\mid 0\leq\beta_i\leq m_i-2\}
$$ 
and 
$\mu=\# I=\Pi_{i=1}^{n+1}(m_i-1)$.
 To calculate the dimensions of the pieces
of $W_\bullet H'$ and $F^\bullet H'$, it is enough to do it for 
$g$. Because in the weighted projective compactification of $\C^{n+1}$
the fibers of $f$ and $g$ are obtained by smooth deformations of each others
and the dimension of the pieces of a mixed Hodge structure is constant under
smooth deformations (see \cite{kk} Chapter 2 \S 3).  We obtain
$$
dim ( Gr^{n+1-k}_FGr_{n+1}^W H')=
\#\{\beta\in I\mid A_\beta=k\}
$$
$$ 
dim(Gr^{n+1-k}_FGr_n^W H')=\#\{\beta\in I \mid  k-1<A_\beta<k\}
$$  
where $A_\beta=\sum _{i=1}^{n+1}
\frac{(\beta_i+1)}{m_i}$. 
Let $P_i$ be the collection of zeros of $\frac{\partial f_i}{\partial x_i}=0$, 
with repetitions according to the  multiplicity, and  $C_i=f_i(A_i)$. 
The set of singularities of $f$ is $P=P_1\times P_2\times\cdots P_{n+1}$ 
and $\sum_{i=1}^{n+1}C_i=\{\sum_{i=1}^{n+1}c_i\mid c_i\in C_i\}$ is the
set of critical values of $f$. The Milnor number of a singularity is
the number of its repetition in $P$.

Before analyzing some examples, let us state a consequence of the Hodge 
conjecture  in the context of this article. We assume that $f$ is a 
strongly tame  polynomial in $\Q[x]$ and $t$ is an algebraic number. 
The Brieskorn module can be redefined over $\Q$ and it turns out that
the Gauss-Manin connection is also defined over $\Q$.   
If the Hodge conjecture is true then a Hodge cycle 
$\delta\in H_n(L_t,\Q)$ satisfies the following property:
For any polynomial differential $n$-form $\omega\in  W_nH'$ (defined over $\Q$) we have
\begin{equation}
\label{9june2005}
\int_\delta\omega\in (2\pi i)^{\frac{n}{2}}\bar \Q
\end{equation}
(See Proposition 1.5 of  Deligne's lecture \cite{dmos}). 
Such a property is proved for 
Abelian varieties of $CM$-type by Deligne.   
Since the main difficulty of the Hodge conjecture lies on construction of 
algebraic cycles, the above statement seems to be much easier to treat 
than the Hodge conjecture itself. 
\begin{exam}\rm 
($f=g=x_1^{m_1}+x_2^{m_2}+\cdots+x_{n+1}^{m_{n+1}},\ $).
 Let
$G:=\Pi_{i=1}^{n+1} G_{m_i}$,
where $G_{m_i}:=\{\epsilon_{m_i}^k \mid k=0,1,\ldots,m_i-1\}$ and  
$\epsilon_{m_i}:=e^\frac{2\pi \sqrt{-1}}{m_i}$ is a primitive root of the unity.
The group $G$ acts on each fiber $L_c$
in the following way:
$$
g: L_c\rightarrow L_c,\ 
(x_1,x_2,\ldots,x_{n+1})\rightarrow
(g_1x_1,g_2x_2,\ldots,g_{n+1}x_{n+1})
$$
where $g=(g_1,g_2,\ldots,g_{n+1})$ is used for
both a vector and a map. Let $I'=\Pi_{i=1}^{n+1}(G_{m_i}-\{1\})$.
We have the one to one map
\begin{equation}
\label{touviajandonaminhacabeca}
I\rightarrow I',\
\alpha \rightarrow 
(\epsilon_{m_1}^{\alpha_1+1},\epsilon_{m_2}^{\alpha_2+1},\ldots,
\epsilon_{m_{n+1}}^{\alpha_{n+1}+1})
\end{equation}
and so we identify $I'$ with $I$ using this map.
Fix a cycle $\delta\in H_n(L_c,\Q)$.
We have
$$
g^*\omega_\beta=g^{\beta+1}\omega_\beta=
\epsilon^{(\alpha+1)(\beta+1)}
\omega_\beta
$$
where 
$$
\epsilon^{(\alpha+1)(\beta+1)}:= 
\epsilon_{m_1}^{(\alpha_1+1)(\beta_1+1)} \epsilon_{m_2}^{(\alpha_2+1)(\beta_2+1)}\cdots \epsilon_{m_{n+1}}^{(\alpha_{n+1}+1)(\beta_{n+1}+1)}
$$
and $g$ corresponds to $\alpha$ by (\ref{touviajandonaminhacabeca}).
We have
\begin{equation}
\label{12.6.03}
\int_{g_*\delta}\omega_\beta=\int_\delta
g^*\omega_\beta=
\epsilon^{(\alpha+1)(\beta+1)}\int_\delta\omega_\beta
\end{equation}
Since $\cup_{\beta\in I} 
\{\delta\in  H_n(L_c,\Q)\mid \int_\delta \omega_\beta=0\}$ does 
not cover $H_n(L_c,\Q)$, we take a cycle $\delta\in H_n(L_c,\Q)$ 
such that $\int _\delta \omega_\beta\neq 0,\ \forall\beta\in I$.
Therefore the period matrix
$P$ in this example is of the form
$E.T$,
where $E= [\epsilon^{(\alpha+1)(\beta+1)}]$ and $T$ is the diagonal matrix 
with $\int_\delta\omega_\beta$ 
in the $\beta\times\beta$ entry. Now the $\omega_\beta,\beta\in I$ form 
a basis of $H'$ and so the period matrix has non zero determinant.
In particular 
the space of Hodge cycles in $H_n(L_c,\Z)$ corresponds to the solutions of
\begin{equation}
\label{9june05}
B.[\epsilon^{(\alpha+1)(\beta+1)}]_{\alpha \in I,\beta\in I_h}=0
\end{equation}
where $B$ is a $1\times\mu$ matrix with integer entries and 
$I_h=\{\beta\in I\mid A_\beta\not\in \N,\ A_\beta<\frac{n}{2}\}$.
This gives an alternative approach for the description of Hodge
cycles for the Fermat variety given by 
Katz, Ogus and Shioda (see \cite{sh79}). 
Note that in their approach one gives an explicit basis of the $\C$-vector 
space generated by Hodge cycles and the elements of such a basis are not 
Hodge cycles. In the description (\ref{9june05}) one can find easily a basis 
of the $\Q$-vector space of Hodge cycles. 
Even if the Hodge conjecture is proved (or disproved), the question of 
constructing an algebraic cycle just with its topological information $B$ 
obtained from (\ref{9june05}) will be  another difficult problem in 
computational algebraic geometry.

\end{exam}
For computations with the next example, we have used 
{\sc Singular} \cite{gps03}. 
\begin{exam}\rm
\label{5.7.03}
($f=x_1^3+x_2^3+\cdots+x_5^3-x_1-x_2$) In this example
$g=x_1^3+x_2^3+\cdots+x_5^3$, $I=\{0,1\}^5$ and $S(t)=27t^3-16t$.
The statement (*) in Remark \ref{21may} is true for 
$$
d_\beta= \left \{ 
\begin{array}{cc}
4 & \beta_1=\beta_2=0 \\
2 & \beta_1=0,\ \beta_2=1  \\
1 & \hbox{ otherwise }
\end{array} \right. .
$$
and so the above data works for all regular values of $f$. 
This follows from the facts that a standard basis of the ideal ${\rm Jacob}(F_t)$, where 
$F_t$ is the homogenization of $f-t$, is given by
$$
2x_1x_0+2x_2x_0+3tx_0^2, x_5^2, x_4^2, x_3^2,
3x_2^2-x_0^2, 3x_1^2-x_0^2, 4x_2x_0^2+3tx_0^3, 
x_0^4
$$
and we have 
{\tiny 
$$
S(t)x_0^4=(\frac{-16}{3}x_2x_0+4tx_0^2) (3x_1^2-x_0^2)+
(\frac{16}{3}x_2x_0+12tx_0^2)(3x_2^2-x_0^2)+
$$
$$
(-8x_1x_2+8x_2^2+6tx_1x_0+6tx_2x_0-9t^2x_0^2)(-2x_1x_0-2x_2x_0+(-3t)x_0^2)
$$
$$
t(4x_2x_0^2+3tx_0^3)=
\frac{4}{9}x_0 (3x_1^2-x_0^2)-
\frac{4}{9}x_0(3x_2^2-x_0^2)+
( \frac{2}{3}x_1-\frac{2}{3}x_2-tx_0)(-2x_1x_0-2x_2x_0-3tx_0^2).
 $$ }
The data $d_\beta=d-1=2$ does not work for the values 
$b=\pm (2/3)^{\frac{3}{2}}$.
In fact for arbitrary $t$ we have 
$x_0x_1x_2=(\frac{9}{8}t^2-\frac{1}{3})x_0^3$ in $V'$, where $V'$ is defined  in Lemma \ref{10:35}.  Now $\nabla^2(\eta)$ is a basis of $Gr^3_FGr_4^W H'_C$. 
We have
{\tiny
$$
\nabla^2(\eta)=\frac{10}{3S(t)^2}( (972t^2-192)x_1x_2\eta+
(-405t^3-48t)x_2\eta
+(-405t^3-48t)x_1\eta
+(243t^4-36t^2+64)\eta)
$$
}which implies the statement in the Introduction. 
%

%

It is remarkable that the
integrals
$\int_\delta \nabla^2(\eta),\ \delta\in H_4(L_t,\Q)$ satisfy the Picard-Fuchs 
equation
\begin{equation}
\label{17Jan}
(27t^3-16t)y''+(81t^2-16)y'+15ty=0
\end{equation}
It is a pull-back of a Gauss hypergeometric equation and so
the integral $\int_\delta \nabla^2(\eta)$ can be expressed in terms
of Gauss hypergeometric series. Since the Hodge conjecture is known 
for cubic hypersurfaces of dimension $4$ by 
C. Clemens, J. P. Murre and S. Zucker (see \cite{zu77}), one can get 
some algebraic relations between the values of such functions on algebraic 
numbers. 
The philosophy of using geometry and obtaining algebraic values of special 
functions  goes back to P. Deligne, 
F. Beukers, J. Wolfart  and many others. 
In  \cite{stho} we have shown that up to a constant, the periods
$\int_\delta\nabla^2\eta,\ \delta\in H_4(L_t,\Q)$ reduce to the periods of
the differential form $\frac{dx}{y}$ on the elliptic curve $E_t: y^2=x^3-3x+z,\ z:=2-\frac{27}{4}t^2$.
It is an interesting observation that $E_b$ associated to 
the fiber  $L_b, \ b=\pm (2/3)^{\frac{3}{2}}$, for which the data 
$d_\beta=d-1$ does not work, is CM. Also the fiber $L_{(2/3)^{\frac{3}{2}}}$ is mapped to 
 $L_{-(2/3)^{\frac{3}{2}}}$ under the automorphism $(x_1,x_2,x_3,x_4,x_5)\mapsto 
 (-x_1,-x_2,\epsilon x_3,\epsilon x_4,\epsilon x_5)$ of the family $f=t$. For the effect of the 
 automorphisms on variation of Hodge structures the reader is referred  to \cite{usu}. 
 We have shown that
the value of the schwarz function $D(0,0,1|z):= -e^{-\pi i \frac{5}{6}} \frac{ F(\frac{5}{6},\frac{1}{6},1|z)}
 {  F(\frac{5}{6},\frac{1}{6},1|1-z)}$  belongs to $\Q(\zeta_3)$ at some $z\in \bar \Q$ if and only if
 $$
 F(\frac{5}{6},\frac{1}{6},1|z) \sim
\frac{1}{\pi^2}\Gamma(\frac{1}{3})^3,\ 
 F(\frac{5}{6},\frac{1}{6},1|1-z)\sim
\frac{1}{\pi^2}\Gamma(\frac{1}{3})^3
$$
where $a\sim b$ means that $\frac{a}{b}\in \bar \Q$.
\end{exam}


Instituto de Matem\'atica Pura e Aplicada, IMPA,\\
 Estrada Dona Castorina, 110,\\
22460-320, Rio de Janeiro, \\
RJ, Brazil.\\
Email: {\tt hossein@impa.br}
\end{document}